\documentclass[letterpaper, 10 pt, conference]{ieeeconf}  % Comment this line out if you need a4paper

\usepackage{amssymb}
\usepackage{amsmath}
\usepackage{graphicx}
\usepackage{multicol}

\IEEEoverridecommandlockouts                              % This command is only needed if
\newcommand{\R}{\mathbb{R}}

\newcommand{\Lp}{\mathbb{L}}

\newcommand{\dis}{\displaystyle}

\newcommand{\vfi}{\varphi}

\newcommand{\sgn}{\mathop{\rm sgn}}

\newcommand{\col}{\mathop{\rm col}}

\newcommand{\be}{\begin{equation}}
\newcommand{\ee}{\end{equation}}
\newcommand{\ba}{\begin{array}}
\newcommand{\ea}{\end{array}}
\newcommand{\RA}{\rightarrow}
\newcommand{\alat}{\mbox{for a.e.}\,t}

\newcommand{\HH}{{\cal H}}

\newtheorem{theorem}{Theorem}
\newtheorem{lemma}{Lemma}
\newtheorem{proposition}{Proposition}

\newtheorem{definition}{Definition}
\newtheorem{hypothesis}{Hypothesis}

\title{\LARGE \bf
Time-Optimal Elevator Control with Higher-Order State Constraints: Analysis and Computation of Boundary Contacts
}

\author{Dmitry~Karamzin$^{1}$, Aleksandra Zhukova$^{2}$, Roman~Chertovskih$^{1}$, A. Pedro Aguiar$^{1}$% <-this % stops a space
\thanks{$^{1}$A.P. Aguiar, R.~Chertovskih, and D.~Karamzin are 
with the Research Center for Systems and Technologies (SYSTEC), ARISE, Faculdade de Engenharia, Universidade do Porto, Rua Dr. Roberto Frias, 4200-465 Porto, Portugal
        {\tt\small pedro.aguiar@fe.up.pt},
        {\tt\small roman@fe.up.pt},
        {\tt\small d.yu.karamzin@gmail.com}}%
\thanks{$^{2}$Aleksandra Zhukova is with the Federal Research Center ``Computer Science and Control'' of Russian Academy of~Sciences, 119333, Vavilova street, 44, Moscow, Russian Federation
        {\tt\small sasha.mymail@gmail.com}}%
}

\begin{document}

\maketitle
\thispagestyle{empty}
\pagestyle{empty}

%%%%%%%%%%%%%%%%%%%%%%%%%%%%%%%%%%%%%%%%%%%%%%%%%%%%%%%%%%%%%%%%%%%%%%%%%%%%%%%%
\begin{abstract}
In this paper, we consider the classical time-optimal elevator problem in the presence of higher-order state constraints. The main contribution is a constructive indirect solution framework based on a Pontryagin Maximum Principle specifically formulated for higher-order constrained systems. First, we derive specialized optimality conditions and analyze the resulting structure of extremal trajectories. Second, we show that the infinite-dimensional optimal control problem can be transformed into a finite-dimensional system of nonlinear algebraic equations involving switching times, boundary-contact times, and multiplier parameters. This provides a computationally efficient procedure for calculating candidate optimal trajectories using standard nonlinear equation solvers. Third, we characterize the geometry of boundary contacts and explain the emergence of contact-chattering phenomena in higher-order constrained systems. In particular, we show that normal extremals cannot evolve along the constraint boundary for a positive amount of time, although sequences of boundary contacts may accumulate indefinitely. The proposed approach is demonstrated on a fourth-order elevator model, and the computed solutions are independently verified using a direct IPOPT-based optimization method.
% In this paper, we consider the elevator problem in its classical formulation and develop a new method to solve it using the Pontryagin maximum principle. The maximum principle is formulated in a special form that accounts for the imposed higher-order state constraints. Its conditions are reduced to a complete system of nonlinear algebraic equations, which is then solved using the standard NLPSolve procedure in Maple. These results are independently confirmed by a direct approach using the IPOPT solver. An interesting feature emerging in this type of control problem is the possibility of theoretically unlimited accumulation of contact with the state-constraint boundary, while the normal trajectory never remains on the boundary for any positive amount of time.

\end{abstract}

%%%%%%%%%%%%%%%%%%%%%%%%%%%%%%%%%%%%%%%%%%%%%%%%%%%%%%%%%%%%%%%%%%%%%%%%%%%%%%%%
\section{INTRODUCTION}\label{Section_1}

Motion planning problems involving bounds on velocity, acceleration, jerk, and higher-order derivatives arise naturally in transportation, robotics, and aerospace applications. In many of these problems, passenger comfort, mechanical limitations, or safety requirements impose constraints not only on the control input but also on the system states and their derivatives. A simple yet representative example is the classical elevator problem, where the objective is to transfer an elevator cabin between two floors in minimum time while satisfying physical and operational limitations. 

Consider a material point moving along a vertical axis. Let the point initially be at rest at some height $x=a$, and suppose that the motion is controlled through its third derivative, $\dddot x=u$, where the control $u$ is bounded. The objective is to transfer the system to a new rest position $x=b$, $a\neq b$, in minimum time. While the unconstrained version of this problem admits a classical solution via the Pontryagin Maximum Principle (PMP), the situation becomes considerably more challenging when state constraints are imposed. In the elevator interpretation, the position is constrained by the shaft boundaries. Additional constraints may naturally be imposed on velocity, acceleration, and higher-order derivatives to ensure safe and comfortable operation. 

Motivated by this example, we consider the following general chain-of-integrators time-optimal control problem on a free time interval $[0,T]$: \be\label{main} \ba{l} \dis T\RA\min,\vspace{0.15cm}\\ \dot x_1(t)=x_2(t),\\ \dot x_2(t)=x_3(t),\\ \quad ...\\ \dot x_{k-1}(t)=x_k(t),\\ \dot x_k(t)=u(t),\\ {\bf x}(0)={\bf x}_0,\\ {\bf x}(T)={\bf x}_T,\\ |u(t)|\le 1\;\alat\in[0,T],\\ x_i(t)\in[h^i_{\min},h^i_{\max}] \;\; \forall\,t\in[0,T],\;\; i=1,\ldots,k, \ea \ee where \[ {\bf x}:=(x_1,x_2,\ldots,x_k)^\top\in\mathbb{R}^k. \] The solution of (\ref{main}) is understood in the sense of a global minimum. Observe that the constraint imposed on $x_i$ is a state constraint of order $k+1-i$ in the terminology of optimal control theory (see, e.g., \cite{Fernando}). Consequently, the constraint on $x_k$ is a standard first-order state constraint, while the constraint on $x_1$ is of order $k$ and requires a substantially more sophisticated analysis. 

Problems of the form (\ref{main}) have attracted considerable attention in the literature. Time-optimal chain-of-integrators systems with full-state constraints were recently studied in \cite{Chinesa}. Related higher-order state-constrained optimal control problems have also been investigated in \cite{Robbins,Dikusar_Milyutin_1989,Zhukova_Karamzin}. A remarkable feature of these problems is the possible occurrence of repeated contacts with the state-constraint boundary. Depending on the boundary conditions, optimal trajectories may exhibit an arbitrarily large number of contact points and, in some cases, an accumulation of contacts, leading to a form of state-constrained chattering. Despite its practical and theoretical relevance, the structure of such solutions and their efficient computation remain challenging issues. 

The contributions of this paper are threefold. First, we specialize the higher-order state-constraint maximum-principle framework of \cite{Fernando} to the time-optimal elevator problem and derive explicit optimality conditions for the corresponding chain-of-integrators model. Second, we establish qualitative properties of the contact set between extremal trajectories and the state-constraint boundary, providing insight into the emergence of finite and accumulating contact sequences and the associated chattering phenomena. Third, we propose a constructive indirect numerical method that transforms the optimality conditions into a complete finite-dimensional system of nonlinear algebraic equations. This formulation enables efficient computation of extremals using standard nonlinear solvers and is independently validated against solutions obtained via direct optimization with IPOPT.

The remainder of the paper is organized as follows. Section~\ref{Section_2} presents the general theoretical framework and the corresponding optimality conditions. Section~\ref{Section_3} analyzes the structure of extremal solutions and boundary-contact phenomena. Section~\ref{Section_4} describes the computational procedure and reports numerical results. Finally, Section~\ref{Section_5} concludes the paper.

\section{GENERAL THEORY}\label{Section_2} Consider the following time-optimal control problem with higher-order state constraints: \be\label{problem} \!\!\!\!\ba{l} \dis T\RA\min,\vspace{0.15cm}\\ \dot x_1(t)= f_1(x_1(t),x_2(t)),\\ \dot x_2(t)= f_2(x_1(t),x_2(t),x_3(t)),\\ \quad...\\ \dot x_{k-1}(t)= f_{k-1}({\bf x}(t)),\\ \dot x_k(t)= f_k({\bf x}(t),u(t)),\\ {\bf p}=({\bf x}(0),{\bf x}(T))\in S,\\ u(t)\in U\;\;\alat\in [0,T],\\ g^j(x_i(t))\le 0\;\forall\, t\in [0,T],\; i=1,\ldots,k,\;j=1,\ldots,l. \ea \ee Here, ${\bf x}=\col(x_1,x_2,\ldots,x_k)$ denotes the state vector in $\R^{n\times k}$, ${\bf p}$ is the endpoint vector in $\R^{2n\times k}$, and $u$ is the control variable taking values in $\R^m$. The sets $S$ and $U$, defining the endpoint constraints and admissible control values, respectively, are assumed to be closed. Furthermore, all mappings appearing in Problem (\ref{problem}) are assumed to be sufficiently smooth. The control function $u(\cdot)$ belongs to $\Lp_\infty([0,T];\R^m)$ and takes values in $U$. A feasible trajectory ${\bf x}(\cdot)$ is an absolutely continuous function satisfying: a) the differential equation $\dot {\bf x}(t)={\bf f}({\bf x}(t),u(t))$ for a.e. $t\in[0,T]$, where ${\bf f}=\col(f_1,f_2,\ldots,f_k)$; b) the endpoint constraints defined by $S$; and c) the inequalities $g^j(x_i(t))\le0$ for all $t\in[0,T]$ and all indices $i,j$, which define state constraints of different orders. A triple $({\bf x}(\cdot),u(\cdot),T)$ satisfying all constraints of Problem (\ref{problem}) is called a feasible process. A feasible process $({\bf \bar x}(\cdot),\bar u(\cdot),\bar T)$ is said to be optimal if $\bar T$ is minimal among all feasible processes. Throughout this section, we assume that Problem (\ref{problem}) admits at least one optimal process $({\bf \bar x}(\cdot),\bar u(\cdot),\bar T)$. Problem (\ref{main}) introduced in the previous section is a particular case of Problem (\ref{problem}). The purpose of this section is to summarize the optimality conditions and auxiliary results that will be used in the subsequent analysis of the elevator problem. For the optimal process, we additionally impose the following regularity assumption. \begin{hypothesis}\label{Hypothesis_1} For every pair of indices $i,j$, the strict inequalities $g^j(\bar x_i(0))<0$ and $g^j(\bar x_i(\bar T))<0$ hold. \end{hypothesis} 

This assumption prevents degeneration of the maximum principle. For a discussion of degeneracy phenomena in optimal control problems with state constraints, see, for example, \cite{Arutyunov_2000} and the references therein. Let $s$ be a positive integer. We denote by $\Gamma_s^{ij}({\bf x},u)$ the $s$-th time derivative of the function $g^j(x_i)$ along trajectories of the differential system in Problem (\ref{problem}), where $s=1,\ldots,k-i+1$. Thus, for every feasible process $({\bf x}(\cdot),u(\cdot),T)$, $$ \Gamma_s^{ij}({\bf x}(t),u(t)) = \frac{d^{\,s}}{dt^s}\Bigl(g^j(x_i(t))\Bigr). $$ For example, when $s=1,2$, $i=1,2$, and $g$ is scalar, one has $\Gamma_1^1=g'f_1$, $\Gamma_1^2=g'f_2$, and $$ \Gamma_2^1=g''[f_1]^2+g'f'_1(f_1,f_2)^*, $$ while $\Gamma_2^2$ is not defined. For larger values of $s$, the explicit expression for $\Gamma_s^{ij}$ becomes increasingly cumbersome, although a recursive representation can still be employed. For convenience, we also define $\Gamma_0^{ij}({\bf x}(t)):=g^j(x_i(t))$. Consider the higher-order Hamilton-Pontryagin function $$ \ba{l} \HH({\bf x},u,\psi,\mu)=\\ \dis\langle\psi,{\bf f}({\bf x},u)\rangle +\sum_{i=1}^{k}\sum_{j=1}^{l} (-1)^i\mu_{ij}\Gamma_{k-i+1}^{ij}({\bf x},u), \ea $$ where $\psi\in(\R^{n\times k})^*$ is the adjoint variable and $\mu=(\mu_{ij})_{i=1,j=1}^{k,l}$ is a matrix-valued multiplier. The optimality conditions for Problem (\ref{problem}) take the following form. \begin{definition}\label{Definition_1} A control process $({\bf \bar x}(\cdot),\bar u(\cdot),\bar T)$ is said to satisfy the maximum principle if there exist a number $\lambda\ge0$, a vector-valued function $\psi\in\mathbb{W}^{\,1}_{\infty}([0,\bar T];(\R^{n\times k})^*)$, and, for each $i=1,\ldots,k$ and $j=1,\ldots,l$, a scalar function $\mu_{ij}\in\mathbb{W}^{k-i}_{\infty}([0,\bar T];\R)$ such that \be\label{t1_1} \dot\psi(t)=-\HH'_{\bf x}(\bar{\bf x}(t),\bar u(t),\psi(t),\mu(t)) \;\;\alat\in[0,\bar T], \ee \be\label{t1_2} \ba{c} \dis\Bigl(\psi(0)+\sum_{j=1}^{l}\sum_{s=1}^{k-i+1} (-1)^s\dot\mu^{(d(s))}(0) (\Gamma^{ij}_{s-1})'(\bar{\bf x}(0)),\\ \dis-\psi(\bar T)\Bigr)\in N_S({\bf \bar p}), \qquad d(s):=k-i+1-s, \ea \ee \be\label{t1_3} \ba{c} \dis\max_{u\in U} \HH(\bar{\bf x}(t),u,\psi(t),\mu(t))\\ = \HH(\bar{\bf x}(t),\bar u(t),\psi(t),\mu(t)) \;\;\alat\in[0,\bar T], \ea \ee \be\label{conservation} \ba{c} \dis\max_{u\in U} \HH(\bar{\bf x}(t),u,\psi(t),\mu(t)) +\\ \dis\sum_{j=1}^{l}\sum_{s=1}^{k-i} (-1)^s\dot\mu^{(d(s))}(t) \Gamma_s^{ij}(\bar{\bf x}(t)) = \lambda \;\;\forall\,t\in[0,\bar T], \ea \ee \be\label{t1_4} \sum_{j=1}^{l}\sum_{i=1}^{k} \int_{0}^{\bar T} g^j(\bar x_i(t)) \,d\dot\mu_{ij}^{(k-i)}(t) = 0, \ee \be\label{t1_5} \lambda+|\psi(0)| + \sum_{j=1}^{l}\sum_{i=1}^{k} \dot\mu_{ij}^{(k-i)}(0) >0. \ee Here, $\mu(t)=(\mu_{ij}(t))_{i=1,j=1}^{k,l}$. Moreover, each function $\dot\mu_{ij}^{(k-i)}$ is decreasing and \be\label{terminal_mu} \mu_{ij}(\bar T) = \dot\mu_{ij}^{(1)}(\bar T) = \cdots = \dot\mu_{ij}^{(k-i)}(\bar T) = 0. \ee \end{definition} \medskip Here, $N_S({\bf p})$ denotes the limiting normal cone to $S$ at ${\bf p}$,~\cite{Mordukhovich_1976}; $\dot\mu^{(r)}$ denotes the derivative of order $r$; and $\mathbb{W}^{\,r}_{\infty}([0,\bar T];\R)$ denotes the space of functions $\vfi$ such that $\dot\vfi^{(r-1)}$ is absolutely continuous and $\dot\vfi^{(r)}$ is bounded whenever $r>0$. For $r=0$, this space coincides with $\Lp_\infty([0,\bar T];\R)$. A control process satisfying the maximum principle is called an extremal. The collection $({\bf \bar x}(\cdot),\bar u(\cdot),\bar T,\lambda,\psi,\mu)$ is referred to as an extended extremal, while $\lambda$, $\psi$, and $\mu$ are the corresponding Lagrange multipliers. When $\lambda>0$, the extremal is called normal. Note that the higher-order derivatives $\dot\mu_{ij}^{(k-i)}$ are well defined at the endpoints due to (\ref{t1_4}) and Hypothesis~\ref{Hypothesis_1}. Indeed, Condition (\ref{t1_4}) implies that $\mu_{ij}(\cdot)$ coincides with a polynomial on every interval $[c,d]$ such that $g^j(\bar x_i(t))<0$ for all $t\in[c,d]$. In particular, $$ \mu_{ij}(t) = \alpha^0_{ij} + \alpha^1_{ij}(t-c) + \cdots + \alpha^{k-i}_{ij}(t-c)^{k-i}, $$ where $\alpha^s_{ij}=\dot\mu^{(s)}(c)/s!$, $s=0,1,\ldots,k-i$. Furthermore, in (\ref{conservation}) we use the fact that $\Gamma_1^{ij}(\bar{\bf x}(t))=0$ whenever $g^j(\bar x_i(t))=0$. These optimality conditions enjoy the following invariance property, which is often useful in practical applications. 

\begin{lemma}\label{Lemma_1} Consider an arbitrary set of Lagrange multipliers $(\lambda,\psi,\mu)$ satisfying the maximum principle. Let $a=(a_1,a_2,\ldots,a_{k-i+1})\in\R^{k-i+1}$ and $c\in[0,1]$. Then, the set of multipliers $$ \Bigl( \lambda, \psi(t) - \sum_{s=1}^{k-i+1} (-1)^s \dot p^{(d(s))}(t) (\Gamma_{s-1}^{ij})'(\bar{\bf x}(t)), \mu_{ij}(t)+p(t) \Bigr), $$ where $p(t)=\sum_{s=0}^{k-i}a_{d(s)}(t-c)^s/s!$ is a polynomial of degree $k-i$, also satisfies conditions (\ref{t1_1}), (\ref{t1_3}), (\ref{conservation}), and (\ref{t1_4}). \end{lemma}

The proof follows the same argument as that presented in \cite{Fernando}. By the invariance principle, choosing $a_s=-\dot\mu_{ij}^{(s-1)}(c)$ yields a new multiplier $\mu_{ij}$ for which, after relabeling, condition (\ref{terminal_mu}) holds with the point $\bar T$ replaced by $c$. The following fundamental result holds. 

\begin{theorem}\label{Theorem_1} Let $(\bar{\bf x}(\cdot),\bar u(\cdot),\bar T)$ be an optimal process in Problem (\ref{problem}). Then, it satisfies the maximum principle. \end{theorem} 

The proof is obtained by reducing Problem (\ref{problem}) to an equivalent problem of twice the state-space dimension and then applying standard optimality conditions to the transformed problem; see \cite{Fernando}. The following two results will be useful in the analysis of the elevator problem.

\begin{theorem}\label{Theorem_2} Let $U$ be convex and compact, $S$ be compact, and let $f_k$ be affine with respect to the control variable $u$ and satisfy a linear growth condition with respect to ${\bf x}$. Then the existence of a feasible process implies the existence of a solution to Problem (\ref{problem}). \end{theorem} 

\begin{theorem}\label{Theorem_3} Let $U$ and $S$ be convex, let ${\bf f}$ be linear with respect to its arguments, and let $g^j$ be convex for all $j=1,\ldots,l$. Suppose there exists a normal extremal $({\bf \bar x}(\cdot),\bar u(\cdot),\bar T,\lambda,\psi,\mu)$. Then, $({\bf \bar x}(\cdot),\bar u(\cdot),\bar T)$ is a solution to Problem (\ref{problem}). \end{theorem} 

Theorem~\ref{Theorem_2} is the classical Filippov existence theorem; see \cite{Filippov_1959}. Theorem~\ref{Theorem_3} is less widely known. However, within the framework of Definition~\ref{Definition_1}, its proof follows from standard convex optimization arguments; see, for example, Chapters~14--16 of \cite{Girsanov}.

\section{ANALYSIS OF OPTIMALITY CONDITIONS IN~THE ELEVATOR PROBLEM}\label{Section_3}

Based on the theory presented in the previous section, let us investigate the solution structure in Problem (\ref{main}). At the same time, for simplicity of the exposition, we assume that $l=1$, and only one scalar state constraint is active, for which $i=s$ for some fixed $s=1,..,k$, for example,
\be\label{state_const}
x_s(t)\le h_{\max}=h^s_{\max}.
\ee
This means that the rest of the constants $h^i_{\min},h^i_{\max}$, $i\ne s$, are considered sufficiently large in the absolute value so the corresponding constraints are nowhere active. The general case of multiple constraints of different orders appears to be more difficult to study, reasoning similar to the one as below can be applied.

Consider the natural assumption that the set of feasible processes is non-empty. Therefore, a solution to problem (\ref{main}) exists by virtue of Theorem \ref{Theorem_2}. This solution, as before, will be denoted by $(\bar{\bf x}(\cdot),\bar u(\cdot),\bar T)$.

Next, for Problem (\ref{main}), all the assumptions of Theorem \ref{Theorem_3} are obviously satisfied. This means that any normal extremal represents a solution. Thus, one can focus on considering only normal extremals when computing and finding solutions to (\ref{main}) using the maximum principle.

Let us apply the optimality conditions stated in the previous section to the elevator problem (\ref{main}).
Denote the extremals in the same way as the solution to the problem, thus, the set $(\bar {\bf x},\bar u,\lambda,\psi,\mu)$ is an extended extremal below. Put $k_s=k-s+1$. Bearing in mind (\ref{state_const}), one has $\Gamma_{k_s}^s({\bf x},u)= u$. Therefore,
$$
\HH({\bf x},u,\psi,\mu) = \sum_{i=1}^{k-1}\psi_i x_{i+1} + \psi_k u + (-1)^{k_s}\mu u,
$$
where $\psi=(\psi_1,\psi_2,...,\psi_k)$.

The optimality conditions take the form
\be\label{f1}
\ba{l}
\dot \psi_1(t) = 0,\\
\dot \psi_2(t) = -\psi_1(t),\\
\dot \psi_3(t) = -\psi_2(t),\\
.....\\
\dot \psi_k(t) = -\psi_{k-1}(t),
\ea
\ee
\be\label{f2}
\ba{c}
\dis\max_{|u|\le 1} \bigl(\psi_k(t) + (-1)^{k_s} \mu(t)\bigr)u\,= \\
\bigl(\psi_k(t) + (-1)^{k_s}\mu(t)\bigr)\bar u(t)\;\;\alat\in [0,\bar T],
\ea
\ee
where $\mu(t)$ satisfies the terminal condition (\ref{terminal_mu}). Moreover, on time intervals where $\bar x_1(t)>0$, this function is a polynomial of degree $k-1$. From (\ref{conservation}) and (\ref{terminal_mu}), it simply follows the equation on $\bar T$,
\be\label{f3}
|\psi_k(\bar T)| = \lambda,
\ee
while, in view of (\ref{f1}), by integrating one has
\be\label{f4}
\psi_k(t) = \sum_{j=0}^{k-1} (-1)^j \alpha_{k-j} \frac{t^j}{j!},
\ee
where $\psi(0)=(\alpha_1,\alpha_2,...,\alpha_k)$.

Consider the closed sets:
$$
\ba{l}
\Sigma:= \bigl\{t\in[0,\bar T]: \bar x_s(t)=h_{\max} \bigr\},\\
\Sigma_0:= \bigl\{t\in\Sigma: \bar x_i(t)=0,\, i=s+1,..,k \bigr\}
\ea
$$
Let us investigate the structure of these sets which is needed for understanding some general properties of a solution. The following simple propositions are valid.

\begin{proposition}\label{Proposition_1}

Let $({\bf \bar x}(\cdot),\bar u(\cdot),\bar T)$ be a normal extremal. Then, the set $\Sigma$ is nowhere dense.

\end{proposition}

Proof. Suppose the contrary, that is, there exist an interval $[a,b]\subseteq (0,\bar T)$ such that $[a,b]\subseteq \Sigma$. Then, obviously, one has $[a,b]\subseteq \Sigma_0$, and $\bar u(t)=0$ on $[a,b]$.
This, however, contradicts the conservation law (\ref{conservation}) when $\lambda>0$. $\quad\square$

Regarding the abnormal extremals, this statement may not already be true.

%However, one can show that it still holds in case of a regular situation, that is, for some open set of the boundary values.

\begin{proposition}\label{Proposition_2}

Let $t_i\in \Sigma$ and $t_i\RA t_0$ for $i\RA\infty$, and let $t_i\ne t_j$ for $i\ne j$. Then, $t_0\in\Sigma_0$.

\end{proposition}

Proof. First, it is obvious that $\bar x_s(t_0)=0$. If we assume that $\bar x_j(t_0)=0$ for $j=s+1,...,l-1$, but $\bar x_l(t_0)\ne 0$ for some $l=s+1,...,k$, then it is easy to demonstrate that $t_0$ is an isolated point in $\Sigma$, which contradicts the assumption made in this proposition. $\quad \square$

Proposition \ref{Proposition_2} and the property of time-optimality imply the following statement.

\begin{proposition}\label{Proposition_3}

Let $s=1$, and $({\bf \bar x}(\cdot),\bar u(\cdot),\bar T)$ be a solution to (\ref{main}). Then,
the set $\Sigma$ is either finite, or it can be represented as the union of a set which has infinitely countable number of isolated points and the singleton set which is the limiting point of the first set from this union. The set $\Sigma_0$ is either empty or a singleton set.

\end{proposition}

For $s>1$, this statement obviously fails to hold as there are simple counter-examples. However, under fairly general assumptions, in case of a regular situation, one can ensure that it is correct when $k-s\ge 2$.

Based on the above, several main cases of interest can be considered. (At the same time, the total number of cases is, of course, greater.)

\begin{itemize}

\item[A)] $\Sigma=[\sigma,\tau]$, where $\sigma<\tau$;

\item[B)] $\Sigma=\{\sigma_1,\sigma_2,...,\sigma_N\}$;

\item[C)] $\Sigma=\{\sigma_i\}_{i=1}^{\infty} \cup \{\sigma\}$, where $\sigma\in\Sigma_0$ is the limiting point.

\end{itemize}

Here, $N$ is a given positive integer number.

Consider Case A). By virtue of Proposition \ref{Proposition_1}, it follows that $\lambda=0$.
Note that $\bar u(t) = 0$ for almost all $t\in \Sigma$, and therefore, the maximum condition implies that
\be\label{psi_k}
\psi_k(t) = (-1)^{k_s+1} \mu(t)\;\;\forall\, t\in [\sigma,\tau].
\ee
At the same time, $\mu(t)=0$ $\forall\,t \in [\tau,\bar T]$, and thereby, the function $\psi_k(t)$ together with its first $k-s-1$ derivatives vanish at the point $\tau$. Thus, in view of (\ref{f3}) and (\ref{f4}), $\psi_k(t)$ can be represented as follows
$$
\psi_k(t) = (t-\tau)^{k-s} p_{s-1}(t),
$$
where $p_{s-1}(t)$ is some polynomial of degree $s-1$ such that $p_{s-1}(\bar T)=0$.

Firstly, consider $s=1$. In this case, $p_0(t)\equiv 0$, and hence, $\psi_k(t)\equiv 0$. Thus, all the multipliers vanish simultaneously which contradict the maximum principle. Then, the case $s=1$ is excluded.

Let $s=2$. Then, $p_1(t)$ is a linear sign-definite function on $[\tau,\bar T]$. This means that $\psi_k(t)$ is sign-definite on this interval. Then, according to the maximum condition, one has that $\bar u(t) = \sgn\psi_k(t)$ for a.e. $t\in [\tau,\bar T]$. However, positive sign implies that $\bar x_s(t)>h_{\max}$ for $t>\tau$ in some neighborhood of $\tau$. This contradicts the arc feasibility. Hence, $\bar u(t)=-1$ for a.e. $t\in [\tau,\bar T]$, and one has
$$
{\bf x}_{T}=\Bigl((-1)^k\frac{\xi^k}{k!},...,\frac{\xi^2}{2},-\xi\Bigr)^*,\;\; \xi=\bar T-\tau.
$$
Thus, the extremals exist only under boundary conditions of the above presented form.

The case $s>2$ is considered similarly, although the formula for ${\bf x}_{T}$ becomes more complicated. At the same time, the conclusion is the same: some specific boundary conditions are needed. We summarize these deductions in the following simple statement.

\begin{proposition}\label{Proposition_4}

For $s=1$, the extremals do not exist in Case A). When $k,s>1$, there may exist only abnormal extremals, but only for a specific set of boundary values which has zero measure.

\end{proposition}

Consider Case B). Consider the polynomial $\psi_k$ whose degree does not exceed $k-1$.
When the degree equals $k-1$, which means that $\alpha_1\ne 0$, one can always try to search for this polynomial in the following decomposed form:
\be\label{roots}
\psi_k(t)= \frac{(-1)^{k-1}\alpha_1}{(k-1)!} (t-r_1)(t-r_2) \,...\, (t-r_{k-1}),
\ee
where $r_1\le r_2\le ...\le r_{k-1}$ are the real roots. Note that $r_i$ may not belong to the interval $[0,\bar T]$.

Let us assume that $k-s>1$, and $s>1$. Firstly, consider the case of a singleton set $\Sigma$, that is, when $N=1$. In this situation, one has $\bar x_s(\sigma_1) = h_{\max}$, and $\bar x_{s+1}(\sigma_1)=0$, as a consequence of Fermat's rule. Let us seek the polynomial $\psi_k$ in the form (\ref{roots}), assuming that all roots $r_i$ are pairwise distinct, and $r_i\in (\sigma_1,\bar T)$, $i=1,2,...,k$. Denote the jump of the function $\dot\mu^{(k-s)}$ at the point $\sigma_1$ by $\delta_1$. Clearly, $\delta_1\le 0$. On $[0,\sigma_1]$, one has $\mu(t) = -\delta_1 (t-\sigma_1)^{k-s}$. Thus, $\mu(t)$ is a parabola, since $k-s\ge 2$. Then, it is simple to see that the polynomial $P(t)=\psi_k(t)+(-1)^{k_s}\mu(t)$ on the interval $[0,\sigma_1]$, depending, of course, on the choice of coefficients $a_i$, can have precisely two roots, which would correspond to switching the control on this interval from $-1$ to $+1$ and back to $-1$, or vice-versa. Indeed, this is so since
\be\label{deg}
\deg\psi_k = k-1> \deg \mu = k-s > 1
\ee
by virtue of the above imposed assumptions.

Thus, there are two independent variables $\delta_1$ and $\sigma_1$, by varying which, within the given limits, one can satisfy two independent equations $\bar x_s(\sigma_1) = h_{\max}$, $\bar x_{s+1}(\sigma_1)=0$. There are also $k$ independent equations which are given by the $k$ boundary conditions $\bar {\bf x}(\bar T)={\bf x}_T$, if $\bar {\bf x}(0)={\bf x}_0$. To solve them, in the normal case, one has $k-1$ independent variables $r_i$, $i=1,...,k-1$ and plus one independent variable $\bar T$. Therefore, the number of equations equals the number of unknowns and thereby, it might be possible to satisfy all the conditions of the maximum principle for a sufficiently large set of boundary values, presumably for some open set of such values. When $\lambda=0$ this is obviously not so as $r_{k-1}=\bar T$, due to (\ref{f3}), and thus, one root is occupied, decreasing the degree of freedom by one.

When $N>1$ the normal extremals can be sought in a similar manner using (\ref{deg}). Let $\delta_i$ denote the jump of $\dot\mu^{(k-s)}$ at $\sigma_i$. Then, there are $2N$ equations given by $\bar x_s(\sigma_i) = h_{\max}$, $\bar x_{s+1}(\sigma_i)=0$, and the same number of unknowns $\sigma_i,\delta_i$, $i=1,...,N$. The roots $r_i$, $i=1,...,k-1$, and the extremal time $\bar T$ are to be found from the given $k$ boundary conditions. Thus, there are $2N+k$ independent variables and the same number of independent algebraic equations. Summing up, by taking into account that for a fixed number $N$ it is possible to consider all combinatorial variants of the roots disposition relative to the contact-points, the following proposition can be formulated.

\begin{proposition}\label{Proposition_5}

In Case B), normal extremals generally exist when $k-s>1$ and $s>1$. Abnormal extremals are possible only for a set of boundary values of measure zero.\footnote{The word `general' here is somewhat vague and requires a clarification. Let $\Theta\subseteq \R^{2n}$ denote the set of all feasible boundary values, that is, the set of all pairs $({\bf x}(0),{\bf x}(T))$ which can be generated by some feasible $u(\cdot)$ and $T$. Then, if solution satisfies Case B), normal extremals exist for almost all points from $\Theta$. Moreover, one can show that the set of boundary values $\Theta_{\rm a}$ for which there exists at least one abnormal extremal is closed, semi-algebraic and has Lebesgue zero measure. This can be shown similar to Case A) as this set can have a direct description in terms of submanifolds (depending on a finite number of arising combinatorial cases). Then, the set $\Theta_{\rm n}=\Theta\setminus \Theta_{\rm a}$ corresponding to the case of exclusively normal extremals is open and has full measure in $\Theta$ provided that $\ell(\Theta)>0$. Here, we also implicitly used that for each element of $\Theta$ there exists a solution to (\ref{main}).}

\end{proposition}

Consider Case C). This case is known as {\it chattering} in literature. It appears rather comprehensive that, in order the function $P(t)$ introduced above to possess infinitely many roots, and thereby, for the control function to have infinitely many switches, one has to assume Condition (\ref{deg}) as a necessary condition. It is necessary indeed since, on each chattering arc, this function must have at least 2 roots to provide two needed switches of the higher-order derivative which are feasible to construct by virtue of the jumps of the $k_s$-order derivative of the multiplier at the contact-points. Thus, the following statement can be proposed.

\begin{proposition}\label{Proposition_6}

Case C) implies Condition (\ref{deg}).

\end{proposition}

Theoretically, nothing prevents the emergence of true infinite chattering. However, the question remains: is this a normal situation or a manifestation of a singularity? Let us formulate the following proposition without proof, that is, simply as a conjecture.

\begin{proposition}\label{Proposition_7}

Case C) implies $\lambda=0$. Moreover, this case is only possible for a specific set of boundary values which has zero measure.

\end{proposition}

Therefore, within a normal situation, only finite chattering is feasible.

\section{COMPUTATIONS}\label{Section_4}

Below, a computational scheme is described to calculate normal extremals in Problem (\ref{main}) by virtue of the  analysis from the previous section. Let us see how, for a given positive integer $N$, to construct an extremal trajectory which contacts $N$ times with the boundary in isolated contact-points. The contact-points are denoted by $\sigma_i$ as before in Case B). Consider the integral of the control dynamical system
$$
I_S^\gamma(t;{\bf s},{\bf b}),\;\;\mbox{where}\;\;{\bf s}=(s_1,s_2,...s_m),
$$
defined on a time interval $S=[t_0,t_1]$ for given $m$ switchings $s_1<s_2<...<s_m$ and starting boundary values ${\bf b}=(b_1,b_2,...,b_k)\in \R^k$. By definition, one has
$$
I_S^\gamma(t;{\bf s},{\bf b})={\bf x}(t),
$$
where the arc ${\bf x}(t)$ satisfies the initial condition ${\bf x}(t_0)={\bf b}$ while the control $u(t)$ is constructed from the switchings $s_i$ as follows:
$$
u(t)=\gamma\cdot(-1)^{i-1},\;\; t\in[s_{i-1},s_i],\;i=1,...,m+1.
$$
Here, $s_0=t_0$, $s_{m+1}=t_1$, $\gamma\in\{-1,+1\}$.

Fix $j\in\{0,1,...,N\}$. Take any $c_j\in(\sigma_j,\sigma_{j+1})$. Here, $\sigma_0=0$, and $\sigma_{N+1}=\bar T$. By virtue of the invariance principle applied at $c_j$, one can assume that Condition (\ref{terminal_mu}) holds at the point $t=c_j$, but not at $t=\bar T$. From the conditions of the maximum principle, one can see that the zeros of
$$
\ba{c}
\Delta(t) = \frac{(-1)^{k-1}\alpha_1}{(k-1)!}\prod_{i=1}^{k-1}(t-r_i)\,-\\
(-1)^{k_s}\sum_{i>j}\delta_i(t-\sigma_i)^{k_s-1}\chi_{[\sigma_i,\bar T]}(t)\,+\\
(-1)^{k_s}\sum_{i<j}\delta_i(t-\sigma_i)^{k_s-1}\chi_{[0,\sigma_i]}(t)
\ea
$$
determine the switching points of the extremal control from $-1$ to $+1$ and back to $-1$, or vice-versa depending on the sign of $\alpha_1$. Here, $r_1,r_2,...,r_{k-1}$, $\delta_i\ge 0$, and $\sigma_i\in (0,\bar T)$, $i=1,...,N$, are independent variables; $\chi_D(t)$ denotes the characteristic function of a subset $D$. Set $\gamma_0={\dis \lim_{t\RA 0+}}\frac{\Delta(t)}{|\Delta(t)|}$.

Consider $r_1<r_2<...<r_{k-1}$ and $r_i\in(\sigma_j,\sigma_{j+1})$ $\forall\, i$ as a baseline scenario.
%Adopt the following local terminology: the switchings $r_i$ are termed explicit switchings, as they are the roots of $\psi_k(t)$.
Suppose that for each $i=1,...,j$ there exist 2 extra switchings $s_{2i-1},s_{2i}\in (\sigma_{i-1},\sigma_i)$ such that $\Delta(s_{2i-1})=\Delta(s_{2i})=0$, while for each $i=j+2,...,N+1$ there exist 2 extra switchings $s_{2i-3},s_{2i-2}\in (\sigma_{i-1},\sigma_i)$ such that $\Delta(s_{2i-3})=\Delta(s_{2i-2})=0$. These switchings are regarded as independent variables. There are also $kN$ variables $a_{ij}= \bar x_j(\sigma_i)$, where $i=1,...,N$ and $j=1,...,k$.
%At the same time, it is clear that $a_{i1}\le 0$ for all $i$.

Thus, one has
$$
\ba{l}
1|_{(\mbox{\scriptsize of}\,\alpha_1)} + k - 1|_{(\mbox{\scriptsize of}\,r_i)} + N|_{(\mbox{\scriptsize of}\,\sigma_i)} + N|_{(\mbox{\scriptsize of}\,\delta_i)} + 2N|_{(\mbox{\scriptsize of}\,s_i)}\,+\\
kN|_{(\mbox{\scriptsize of}\,a_{ij})} + 1|_{(\mbox{\scriptsize of}\,\bar T)} = k+(k+4)N+1.
\ea
$$
independent variables. Now, it is needed to find the same number of equations.
The relation
$$
I^{\gamma_0}_{[0,\sigma_1]}(\sigma_1;(s_1,s_2),{\bf x}_0) = {\bf a}_1,
$$
produces the first $k$ equations. Here, ${\bf a}_i:=(a_{i1},a_{i2},...,a_{ik})$.
 %$D_s$ is the matrix which truncates the first $s-1$ elements of a vector.
Then, for each $i=1,2,...,j-1$, one has the $k$ equations
$$
I^{\gamma_0}_{[\sigma_i,\sigma_{i+1}]}(\sigma_{i+1};(s_{2i+1},s_{2i+2}),{\bf a}_i) = {\bf a}_{i+1}.
$$
When $i=j$, the $k$ equations are as follows:
$$
I^{\gamma_0}_{[\sigma_j,\sigma_{j+1}]}(\sigma_{j+1};(r_1,r_2,...,r_{k-1}),{\bf a}_j) = {\bf a}_{j+1}.
$$
For $i=j+1,j+2,...,N-1$, the $k$ equations are
$$
I^{\gamma_1}_{[\sigma_i,\sigma_{i+1}]}(\sigma_{i+1};(s_{2i-1},s_{2i}),{\bf a}_i) = {\bf a}_{i+1},
$$
where $\gamma_1=\gamma_0\cdot(-1)^{k-1}$. Finally, for $i=N$, these equations are
$$
I^{\gamma_1}_{[\sigma_N,\bar T]}(\bar T;(s_{2N-1},s_{2N}),{\bf a}_N) = {\bf x}_T.
$$
This already gives $k + kN$ equations. The equalities $\Delta(s_i)=0$ additionally yield $2N$ equations. More $2N$ equations are given by the contact condition: $a_{is}=h_{\max}^s$, $a_{i(s+1)}=0$. Thus, the total number of equations equals $k+N(k+4)$, and thereby it remains to find just one more independent equation. This equation is given by formula (\ref{f3}) which, under the variable change of Lemma \ref{Lemma_1} with $c=c_j$, takes the following form:
$$
\lambda = \frac{|\alpha_1|}{(k-1)!}\prod_{i=1}^{k-1}|c_j-r_i|.
$$
One can consider $\lambda=1$. Then, we arrive at the ultimate equation and therefore, the total number of unknowns equals the total number of equations. This means that one can try solving numerically these $k+N(k+4)+1$ algebraic equations with respect to the available variables listed above. At the same time, the number $j$ can vary from 0 to $N$ which means that for each $j$ the same numerical procedure can be applied increasing the searching range. Note that when $j=0$, one can take $c=0$, while for $j=N$, $c=\bar T$. Note that in the last equation, there are always 2 roots for $\alpha_1$.\footnote{Also note that in the case of symmetric boundary values, that is, when the arc and its derivatives enjoy the property of central or axial symmetry, the presented scheme collapses for the case of odd $N$. Then, one needs to join the two central intervals putting all $r_i$ therein.}

\begin{figure}[h]

 \centering

 \includegraphics[width=70mm]{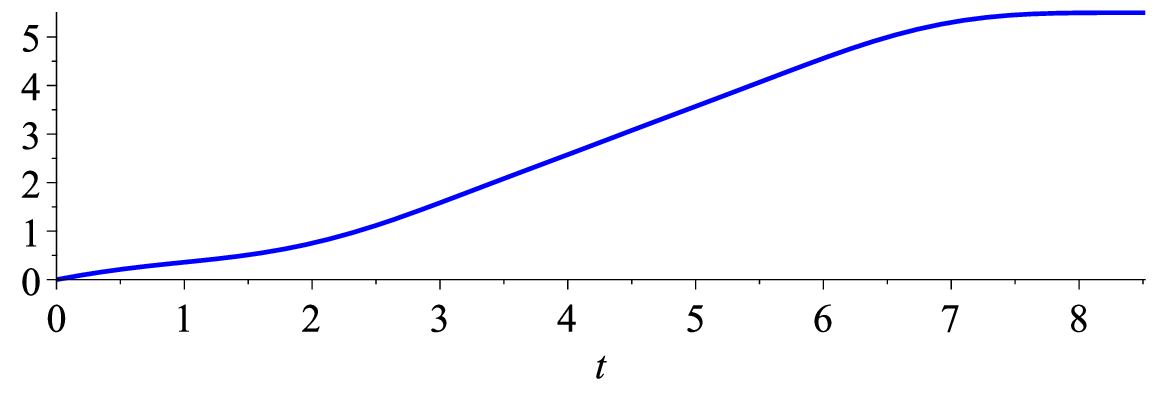}

 \caption{The graph of $\bar x_1(t)$.}

 \label{figure_1}

 \end{figure}

\begin{figure}[h]

 \centering

 \includegraphics[width=70mm]{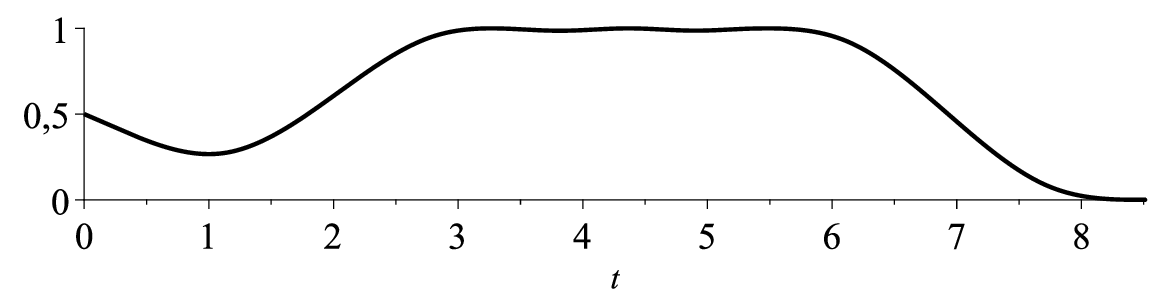}

 \caption{The graph of $\bar x_2(t)$.}

 \label{figure_2}

 \end{figure}

\begin{figure}[h]

 \centering

 \includegraphics[width=70mm]{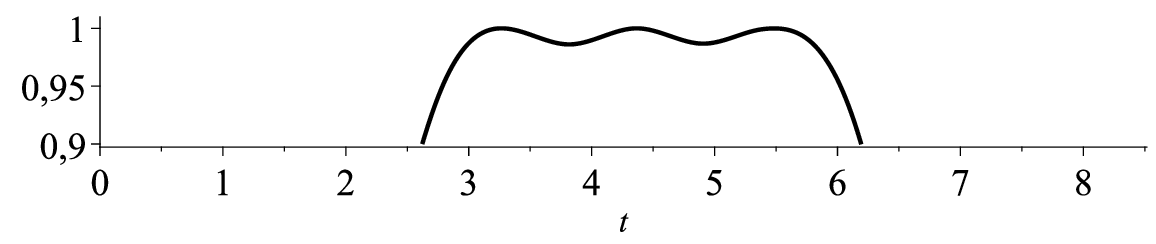}

 \caption{Zoom in of $\bar x_2(t)$.}

 \label{figure_3}

 \end{figure}

\begin{figure}[h]

 \centering

 \includegraphics[width=70mm]{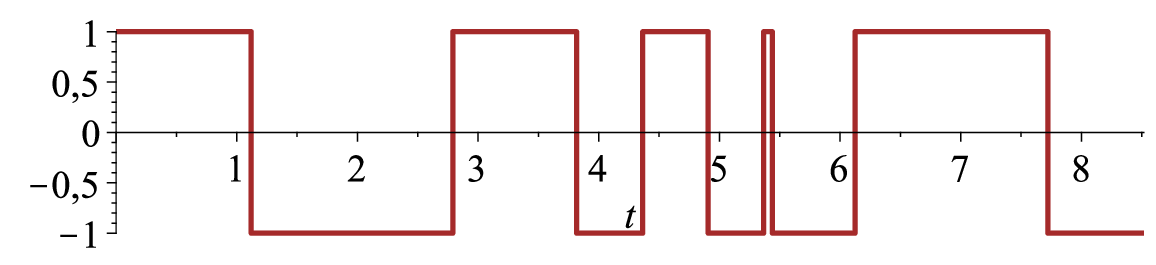}

 \caption{The graph of control $\bar u(t)$.}

 \label{figure_4}

 \end{figure}

Below as an example, the above computational scheme is applied to calculate the solution to the elevator problem in case $k=4$, $s=2$, $N=3$. In Figure \ref{figure_1} \ref{figure_2}, the optimal height-position and velocity are shown when the boundary values are ${\bf \bar x}(0) = (0,0.5,-0.3,-0.2)$, ${\bf \bar x}(T) = (5.5,0,0,0)$, while $h^2_{\max}=1$. The accuracy of solving equations (FuncTol) is $\approx 10^{-9}$. In Figure \ref{figure_3}, one can observe three contact-points. In Figure \ref{figure_4}, the optimal control is shown.

%\begin{figure}[h]

% \centering

%\begin{multicols}{2}

% \includegraphics[width=45mm]{PICS/Pic_1.eps}
% \includegraphics[width=45mm]{PICS/Pic_2.eps}

% \end{multicols}

% \caption{The graphs of $\bar x_1(t)$ and $\bar x_2(t)$.}

% \label{figure_1}

% \end{figure}

The calculated values are given in the next table.
$$
{\footnotesize
\ba{c}
s1 = 1.119704487,\\
s2 = 2.790212626,\\
s3 = 3.816151814,\\
s4 = 4.364275685,\\
s5 = 4.906826920,\\
s6 = 5.366863762,\\
r1 = 5.439041929,\\
r2 = 6.124079611,\\
r3 = 7.722170744,
\ea
\qquad
\ba{c}                        
\alpha_1 = 1.025058536,\\
\delta_1 = 0.859450136,\\
\delta_2 = 0.723222015,\\
\delta_3 = 0.891043786,\\
\sigma_1 = 3.265880736,\\
\sigma_2 = 4.366439917,\\
\sigma_3 = 5.533295525,\\
\mbox{and}\\
\bar T = 8.516928415.
\ea
}
$$
These computations have been independently confirmed using IPOPT-solver and a similar solution-structure found. The IPOPT-optimal time is about $0.01$ less than $\bar T$ shown above.

\section{CONCLUSIONS}\label{Section_5}

This paper presents the general theoretical framework for the analysis of the time-optimal elevator control problem and the corresponding optimality conditions. The presented analysis of the structure of extremal solutions and boundary-contact phenomena suggests the computational procedure to obtain the optimal control and trajectories. The computational scheme has been applied to a case of 4-dimensional elevator problem with constrained velocity and the resulting solutions have been visualized on graphs.

\addtolength{\textheight}{-12cm}   % This command serves to balance the column lengths

%%%%%%%%%%%%%%%%%%%%%%%%%%%%%%%%%%%%%%%%%%%%%%%%%%%%%%%%%%%%%%%%%%%%%%%%%%%%%%%%

%%%%%%%%%%%%%%%%%%%%%%%%%%%%%%%%%%%%%%%%%%%%%%%%%%%%%%%%%%%%%%%%%%%%%%%%%%%%%%%%

%%%%%%%%%%%%%%%%%%%%%%%%%%%%%%%%%%%%%%%%%%%%%%%%%%%%%%%%%%%%%%%%%%%%%%%%%%%%%%%%
% \section*{ACKNOWLEDGMENT}

% The authors acknowledge the support from R\&D Unit SYSTEC -- UID/EEA/00147 -- POCI-01-0145-FEDER-006933 funded by ERDF $|$ COMPETE2020 $|$ FCT/MEC $|$ PT2020, and also from the STRIDE project NORTE-01-0145-FEDER-000033, by ERDF $|$ NORTE 2020, Porto University, Portugal.

%%%%%%%%%%%%%%%%%%%%%%%%%%%%%%%%%%%%%%%%%%%%%%%%%%%%%%%%%%%%%%%%%%%%%%%%%%%%%%%%


\begin{thebibliography}{99}

\bibitem{Fernando}
D. Karamzin, F.L. Pereira, On higher-order state constraints, SIAM J. Control Optim., vol.~61, no.~4, 2023, pp.~1913--1933.

\bibitem{Chinesa}
Y. Wang, C. Hu, Z. Li, S. Lin, S. He, Y. Zhu, Time-Optimal Control for High-Order Chain-of-Integrators Systems With Full State Constraints and Arbitrary Terminal States, IEEE Transactions on Automatic Control, vol. 70, no. 3, 2025, pp. 1499--1514.

\bibitem{Robbins}
H. Robbins, Junction phenomena for optimal control with state-variable inequality constraints of third order,
J. Optim. Theory Appl., vol. 31, 1980, pp. 85--99.

\bibitem{Dikusar_Milyutin_1989}
V.V. Dikusar, A.A. Milyutin. Qualitative and numerical methods in the maximum principle. Moscow, Nauka, 1989.

\bibitem{Zhukova_Karamzin}
A.A. Zhukova, D.Y. Karamzin, Calculation of Extremals in an Optimal Control Problem with a Higher-Order State Constraint, Comput. Math. and Math. Phys., vol. 65, 2025, pp. 2838--2853.

\bibitem{Arutyunov_2000}
A.V. Arutyunov, Optimality conditions, Abnormal and Degenerate Problems. Mathematics and Its Application, Kluwer AP, 2000.

\bibitem{Mordukhovich_1976}
B.S. Mordukhovich, Maximum principle in the problem of time optimal response with nonsmooth constraints, J. Appl. Math. Mech., vol. 40, no. 6, 1976, pp. 960--969.

\bibitem{Filippov_1959}
A.F. Filippov, On certain problems of optimal regulation, Bull. of Moscow State University, ser. Math. and Mech., 1959, pp. 25--38.

\bibitem{Girsanov}
I.V. Girsanov, Lectures on the Mathematical Theory of Extremal Problems. Moscow University Press, 1970.

\bibitem{IPOPT}
A. W\"{a}chter, L.T. Biegler, On the implementation of an interior-point filter line-search algorithm for large-scale nonlinear programming, Math. Program., vol. 106, no. 1, 2006, pp. 25--57.

\end{thebibliography}
\end{document}